# RATES OF CONTRACTION OF POSTERIOR DISTRIBUTIONS BASED ON GAUSSIAN PROCESS PRIORS

By A. W. van der Vaart and J. H. van Zanten[1]

*Vrije Universiteit Amsterdam*

We derive rates of contraction of posterior distributions on nonparametric or semiparametric models based on Gaussian processes. The rate of contraction is shown to depend on the position of the true parameter relative to the reproducing kernel Hilbert space of the Gaussian process and the small ball probabilities of the Gaussian process. We determine these quantities for a range of examples of Gaussian priors and in several statistical settings. For instance, we consider the rate of contraction of the posterior distribution based on sampling from a smooth density model when the prior models the log density as a (fractionally integrated) Brownian motion. We also consider regression with Gaussian errors and smooth classification under a logistic or probit link function combined with various priors.

**1. Introduction.** Gaussian processes have been adopted as building blocks for constructing prior distributions on infinite-dimensional statistical models in several settings. For instance, in the setting of nonparametric density estimation, a prior distribution on a collection of probability densities (relative to a measure $\nu$) can be defined structurally as the random density

$$(1.1) \qquad \frac{e^{W_x}}{\int e^{W_x} \, d\nu(x)},$$

where $(W_x : x \in \mathcal{X})$ is a Gaussian process indexed by the sample space $\mathcal{X}$ of the observations. The Gaussian process is exponentiated to force the prior to charge only nonnegative functions, and is next renormalized to integrate to unity. Several other constructions have also been considered in the literature, for density estimation as well as other statistical problems; see Section 3 and [4, 5, 9, 12, 15, 19, 20, 21, 22, 26, 29, 34]. The book [27] makes a connection to

Received May 2006; revised March 2007.

[1]Supported by the Netherlands Organization for Scientific Research NWO.

*AMS 2000 subject classifications.* 60G15, 62G05.

*Key words and phrases.* Rate of convergence, Bayesian inference, nonparametric density estimation, nonparametric regression, classification.

machine learning and the website <http://www.gaussianprocess.org> lists additional references.

Given a prior and observations, Bayes' rule yields a posterior distribution on the parameter space. In the frequentist set-up, in which the data are sampled from a fixed "true" distribution and the amount of information in the data increases indefinitely, the corresponding posterior distributions often contract to the fixed true distribution, which is referred to as *posterior consistency*. In this paper, we study the *rate of contraction* of the posterior distribution relative to global metrics on the parameters.

In most cases, the Gaussian process can be viewed as a tight Borel measurable map in a Banach space, for instance, a space of continuous functions or an $L_p$-space. It is well known that the support of a centered (i.e., zero-mean) version of such a process (the smallest closed set having probability one under the induced measure) is equal to the closure of the *reproducing kernel Hilbert space* (RKHS) of the covariance kernel of the process. Because the posterior distribution necessarily puts all of its mass on the support of the prior, it follows that consistency can be valid only if the parameter $w_0$ defining the true distribution of the data belongs to this support. In the present paper, we prove that the rate of contraction in that case is expressible in terms of the function

$$(1.2) \qquad \phi_{w_0}(\varepsilon) = \inf_{h \in \mathbb{H}: \|h - w_0\| < \varepsilon} \|h\|_{\mathbb{H}}^2 - \log \Pr(\|W\| < \varepsilon).$$

In this expression, $\|\cdot\|$ is the norm of the Banach space in which the Gaussian process $W$ takes its values, $\mathbb{H}$ is the reproducing kernel Hilbert space of the process and $\|\cdot\|_{\mathbb{H}}$ the RKHS-norm. If the norm $\|\cdot\|$ on the sample space of the process "combines correctly" with the norm on the parameter space and $n$ describes the informativeness of the data in the usual way, then the posterior contracts at the rate $\varepsilon_n \to 0$, satisfying

$$(1.3) \qquad \phi_{w_0}(\varepsilon_n) \le n\varepsilon_n^2.$$

This is the case, for instance, in density estimation on the unit interval with the log Gaussian process prior given in (1.1) and $\|\cdot\|$ the uniform norm given by $\|w\| = \sup\{|w(x)| : x \in \mathcal{X}\}$. This is also the case for regression and classification, with appropriate norms, as shown below. The rate of contraction $\varepsilon_n$ thus depends on the position of the true parameter $w_0$ relative to the RKHS and the amount of mass $\Pr(\|W\| < \varepsilon)$ that the prior distribution puts in small balls around zero. In Section 4, we compute these quantities for a range of priors.

For instance, we prove that log Gaussian densities given in (1.1), combined with Brownian motion, yield a rate of contraction of $n^{-1/4}$ whenever the logarithm of the true density is $\alpha$-smooth for some $\alpha \ge 1/2$ and yield the slower rate $n^{-\alpha/2}$ for $0 < \alpha \le 1/2$. That higher smoothness ($\alpha$ large) does not improve the rate of contraction is disappointing, but perhaps not surprising, given that the Brownian paths themselves are 1/2-smooth: the data are not capable of smoothing out the prior roughness of the sample paths. Other, more smooth, Gaussian priors give better rates for smooth truths (depending on their RKHS), but worse rates for rough truths. In Section 4, we exhibit, for every possible smoothness level $\alpha$, Gaussian priors that give the optimal rate of contraction if the true parameter possesses regularity $\alpha$.

The function $\phi_{w_0}$ displayed in (1.2) may seem complicated at first. In fact, it can be handled for many examples. In particular, the probability $\Pr(\|W\| < \varepsilon) = e^{-\phi_0(\varepsilon)}$ is known as the *small ball probability* for $\varepsilon \downarrow 0$ and has been studied in many papers in the probability literature (see [24] or the extensive bibliography compiled by M. Lifshits on http://www.proba.jussieu.fr/pageperso/smalldev/biblio.pdf). In Section 4, we discuss a number of examples. The centered small ball probability exponent $\phi_0(\varepsilon)$ puts a limit on the rate of contraction that depends only on the prior, while the decentered small ball probability shows how this rate might deteriorate by the positioning of the true parameter $w_0$ relative to the support of the prior.

The paper is organized as follows. In Section 2, we recall the definition of the RKHS of a Gaussian process and state theorems on the concentration of Gaussian processes that are the basis of the remainder of the paper. In Section 3, we state our main results on posterior concentration for a number of statistical settings. Next, in Section 4, we discuss a number of special Gaussian processes and derive the rates of posterior contraction corresponding to true parameters of various regularity levels. Section 5 contains the proofs.

The notation $\lesssim$ is used for "smaller than or equal to a universal constant times" and $\asymp$ means "proportional up to constants." We let $\|\cdot\|_{p,\nu}$ denote the norm of $L_p(\nu)$, the space of measurable functions with $\nu$-integrable $p$th absolute power. Furthermore, $h(f,g) = \|\sqrt{f} - \sqrt{g}\|_{2,\nu}$ is the Hellinger distance, $K(f,g)$ the Kullback–Leibler divergence, and $V(f,g) = \int (\log(f/g))^2 f \, d\nu$. If the dominating measure $\nu$ is Lebesgue measure, then it may be omitted in the notation. The notation $C[0,1]$ is used for the space of continuous functions $f : [0,1] \to \mathbb{R}$ endowed with the uniform norm and, for $\beta > 0$, we let $C^\beta[0,1]$ denote the Hölder space of order $\beta$, consisting of the functions $f \in C[0,1]$ that have $\underline{\beta}$ continuous derivatives for $\underline{\beta}$ the biggest integer strictly smaller than $\beta$ with the $\underline{\beta}$th derivative $f^{(\underline{\beta})}$ being Lipschitz continuous of order $\beta - \underline{\beta}$. Finally, $H^k[0,1]$ denotes the Sobolev space of functions $f : [0,1] \to \mathbb{R}$ that are $k - 1$ times continuously differentiable with absolutely continuous $(k-1)$th derivative that is the integral of a function $f^{(k)} \in L_2[0,1]$ and $\ell^\infty(\mathcal{X})$ is the space of bounded functions $z : \mathcal{X} \to \mathbb{R}$ with the uniform norm $\|z\|_{\mathcal{X}} = \sup\{|z(x)| : x \in \mathcal{X}\}$ (also written as $\|z\|_\infty$).

**2. Gaussian priors.** In this section, we first recall the definition of the RKHS and next formulate results on the support of Gaussian processes which will be basic to the results on rates of posterior contraction in the next section. The proofs of the results in this section are deferred to Section 5. Relevant results on RKHS are scattered throughout the literature. Van der Vaart and Van Zanten [30] reviews facts that are relevant to the present applications.

The definition of an RKHS that is most appropriate for the results in this paper concerns Gaussian random elements seen as Borel measurable maps in a Banach space. A Borel measurable random element $W$ with values in a separable Banach space $(\mathbb{B}, \|\cdot\|)$ is called *Gaussian* if the random variable $b^*W$ is normally distributed for any element $b^*$ of the dual space $\mathbb{B}^*$ of $\mathbb{B}$ and it is called *zero-mean* if the mean of every such variable $b^*W$ is zero. The *reproducing kernel Hilbert space* (*RKHS*) $\mathbb{H}$ attached to $W$ is the completion of the range $S\mathbb{B}^*$ of the map $S: \mathbb{B}^* \to \mathbb{B}$ defined by

$$Sb^* = \mathrm{E}Wb^*(W), \qquad b^* \in \mathbb{B}^*,$$

for the inner product

$$\langle Sb_1^*, Sb_2^* \rangle_{\mathbb{H}} = \mathrm{E}b_1^*(W)b_2^*(W).$$

The element $Sb^* \in \mathbb{B}$ is the Pettis integral of the $\mathbb{B}$-valued random element $Wb^*(W)$—an element of $\mathbb{B}$ such that $b_2^*(Sb^*) = \mathrm{E}b_2^*(W)b^*(W)$ for every $b_2^* \in \mathbb{B}^*$ (cf. [18], page 42). It can be shown the RKHS-norm on the set $S\mathbb{B}^*$ is stronger than the original norm, so the RKHS $\mathbb{H}$, the completion of the set $S\mathbb{B}^*$ under the RKHS-norm, can be identified with a subset of $\mathbb{B}$.

A zero-mean Gaussian stochastic process $W = (W_t : t \in T)$ defined on some probability space $(\Omega, \mathcal{U}, \mathrm{Pr})$ and indexed by an arbitrary set $T$ with bounded sample paths $t \mapsto W_t$ can be viewed as a map into the Banach space $\ell^\infty(T)$. If it is Borel measurable and has separable range, then its RKHS is defined above. It can be shown (e.g., [30]) that this RKHS can be identified with the completion of the set of maps

$$(2.1) \qquad t \mapsto \sum_i \alpha_i K(s_i, t) = \mathrm{E}W_t H, \qquad H = \sum_i \alpha_i W_{s_i},$$

under the inner product

$$\langle \mathrm{E}W.H_1, \mathrm{E}W.H_2 \rangle_{\mathbb{H}} = \mathrm{E}H_1 H_2.$$

Here, $K(s,t) = \mathrm{E}W_sW_t$ is the covariance function of the process and $H$ ranges over all finite linear combinations. This completion is precisely the set of functions $t \mapsto \mathrm{E}W_t H$ with $H$ ranging over the closure of the set of linear combinations $H = \sum_i \alpha_i W_{s_i}$ in $L_2(\Omega, \mathcal{U}, \mathrm{Pr})$.

For $\varepsilon > 0$, let $N(\varepsilon, B, d)$ denote the minimum number of balls of radius $\varepsilon$ needed to cover a subset $B$ of a metric space with metric $d$ (cf. [31]).

THEOREM 2.1. *Let $W$ be a Borel measurable, zero-mean Gaussian random element in a separable Banach space $(\mathbb{B}, \|\cdot\|)$ with RKHS $(\mathbb{H}, \|\cdot\|_{\mathbb{H}})$ and let $w_0$ be contained in the closure of $\mathbb{H}$ in $\mathbb{B}$. For any numbers $\varepsilon_n > 0$ satisfying* (1.3) *for $\phi_{w_0}$ given by* (1.2), *and any $C > 1$ with $e^{-Cn\varepsilon_n^2} < 1/2$, there exists a measurable set $B_n \subset \mathbb{B}$ such that*

$$\text{(2.2)} \qquad \log N(3\varepsilon_n, B_n, \|\cdot\|) \leq 6Cn\varepsilon_n^2,$$

$$\text{(2.3)} \qquad \Pr(W \notin B_n) \leq e^{-Cn\varepsilon_n^2},$$

$$\text{(2.4)} \qquad \Pr(\|W - w_0\| < 2\varepsilon_n) \geq e^{-n\varepsilon_n^2}.$$

The three assertions of this theorem can be matched one-to-one with the assumptions of general results on rates of posterior contraction (e.g., Theorem 2.1 of [8]), except that the assertions here use the norm of the Banach space, whereas the conditions for the posterior rates are in terms of metrics or discrepancies appropriate to the statistical problem under consideration. The rate of contraction $\varepsilon_n$ is obtained as soon as the latter metrics are comparable to the norm. This is shown to be the case for various statistical settings in the next section.

The preceding theorem is meant to be used as an asymptotic result as $n \to \infty$, but is, in fact, a statement for every fixed $n$. The Gaussian process $W$ and the "true" parameter $w_0$ may therefore also be taken to be dependent on $n$, as long as the corresponding RKHS and function $\phi_{w_n}$ are also taken to be dependent on $n$.

In the context of sequences of Gaussian processes that approximate a fixed process, such as truncated Fourier series, working with a sequence of concentration functions would be unnecessarily cumbersome. We have the following refinement, which shows that we can use the concentration function of the limit process in such cases.

THEOREM 2.2. *Let $W^n$ be Borel measurable, zero-mean, jointly Gaussian random elements in a separable Banach space $(\mathbb{B}, \|\cdot\|)$ such that $10\mathrm{E}\|W^n - W\|^2 \leq 1/n$ for a Gaussian process $W$. Let $(\mathbb{H}, \|\cdot\|_{\mathbb{H}})$ be the RKHS of $W$ and assume that $w_0$ is contained in the closure of $\mathbb{H}$ in $\mathbb{B}$. For any numbers $\varepsilon_n > 0$ satisfying $n\varepsilon_n^2 \geq 4\log 4$ and* (1.3) *with $\phi_{w_0}$ given by* (1.2), *and any $C > 4$ with $e^{-Cn\varepsilon_n^2} < 1/2$, there exists a measurable set $B_n \subset \mathbb{B}$ such that* (2.2), (2.3) *and* (2.4) *hold with $W$ replaced by $W^n$ and $\varepsilon_n$ replaced by $2\varepsilon_n$.*

The sum $W = \sum_i W^i$ of finitely many independent Gaussian processes $W^i$ is itself a Gaussian process. It appears that it is not always easy to obtain its RKHS from the RKHS's of the components $W^i$. However, the concentration function of $W$ can easily be obtained from the concentration functions of the components.

THEOREM 2.3. *Let $W = \sum_{i \in I} W^i$ be the sum of finitely many independent Borel measurable, zero-mean Gaussian random elements in a separable Banach space $(\mathbb{B}, \|\cdot\|)$ with concentration functions $\phi_{w^i}^i$ for given $w^i \in \mathbb{B}$. Then, the concentration function $\phi_w$ of $W$ around $w = \sum_{i \in I} w^i$ satisfies*

$$\phi_w(\varepsilon|I|) \le 2 \sum_{i \in I} \phi_{w^i}^i(\varepsilon/2).$$

The theorem applies, in particular, to a sum $V + W$ where $W$ possesses the desired properties (2.2), (2.3) and (2.4) and $V$ is more concentrated at zero than $W$, in the sense that $\Pr(\|V\| \le \varepsilon) \ge \Pr(\|W\| \le \varepsilon)$ for every $\varepsilon > 0$. The theorem with $W^1 = V$, $w^1 = 0$ and $W^2 = W$ shows that $V$ will not destroy good properties of $W$ in that case.

It is natural to scale a Gaussian process so that its fluctuations are of the same order of magnitude as the fluctuations thought to exist in the true parameter $w_0$. Lacking sufficient prior insight regarding $w_0$, one might use a hierarchical prior of the form $AW$, where the scale parameter $A$ is chosen from some distribution on $(0, \infty)$, independent of the Gaussian process $W$. The preceding results extend to this prior if the support of the prior for $A$ is bounded above. (The rate deteriorates if the scale parameter is not bounded away from infinity. We do not discuss this case here.)

THEOREM 2.4. *Let $W$ be a Borel measurable, zero-mean Gaussian random element in a separable Banach space $(\mathbb{B}, \|\cdot\|)$ independent of the random variable $A$ that takes its values in an interval $(0, K] \subset (0, \infty)$. Let $w_0$ be contained in the closure of the RKHS $\mathbb{H}$ of $W$ in $\mathbb{B}$. Let $k < 1 < K$. For any numbers $\varepsilon_n > 0$ satisfying (1.3) for $\phi_{w_0}$ given by (1.2), and any $C > 1$ with $e^{-Cn\varepsilon_n^2} < 1/2$, there exists a measurable set $B_n \subset \mathbb{B}$ such that*

$$(2.5) \qquad \log N(3K\varepsilon_n, B_n, \|\cdot\|) \le 6Cn\varepsilon_n^2,$$

$$(2.6) \qquad \Pr(AW \notin B_n) \le e^{-Cn\varepsilon_n^2},$$

$$(2.7) \qquad \Pr(\|AW - w_0\| < 2K\varepsilon_n) \ge \Pr(A \ge k)e^{-n\varepsilon_n^2/k^2}.$$

**3. Main results on posterior contraction.** Gaussian processes can be used as building blocks for constructing priors on function spaces in various ways and in several statistical settings. In order for our general approach to apply, appropriate metrics on the set of distributions of the observations must correspond to the norm of the Banach space in which the Gaussian process takes its values. In this section, we describe several cases where this desirable situation is achieved. These are motivated by implementations in the literature and do not form an exhaustive set.

3.1. *Density estimation.* Suppose that we observe an i.i.d. sample $X_1, \ldots,$ $X_n$ from a density $p_0$ relative to a measure $\nu$ on a measurable space $(\mathcal{X}, \mathcal{A})$. Consider a prior distribution on the set of $\nu$-densities defined structurally as $p_W$ for a Gaussian process $W$ and, for $p_w$, the function defined by

$$p_w(x) = \frac{e^{w_x}}{\int_{\mathcal{X}} e^{w_y} \, d\nu(y)}.$$

(The notation $p_0$ now denotes both the true density and the density $p_w$ with $w = 0$.) Implementations of this prior were considered in [19, 20, 22] addressing, for instance, the computation of the posterior mean.

Assume that $W$ has bounded sample paths and can be viewed as a Borel measurable map in the space $\ell^\infty(\mathcal{X})$ of bounded functions $z : \mathcal{X} \to \mathbb{R}$ equipped with the uniform norm. The following theorem shows that the rate of contraction for log Gaussian prior densities is determined exactly as in (1.2)–(1.3), with $w_0 = \log p_0$.

THEOREM 3.1. *Let $W$ be a Borel measurable, zero-mean, tight Gaussian random element in $\ell^\infty(\mathcal{X})$. Suppose that $w_0 = \log p_0$ is contained in the support of $W$ and let $\phi_{w_0}$ be the function in (1.2) with $\| \cdot \|$ the uniform norm on $\ell^\infty(\mathcal{X})$. Then, the posterior distribution relative to the prior $p_W$ satisfies $\mathrm{E}_0 \Pi_n(p_w : h(p_w, p_0) > M \varepsilon_n | X_1, \ldots, X_n) \to 0$ for any sufficiently large constant $M$ and $\varepsilon_n$ given by (1.3).*

PROOF. The proof of the theorem is based on Theorem 2.1 of [8] and a comparison of the Hellinger and Kullback–Leibler distances between log Gaussian prior densities to the uniform distance on the Gaussian process, as in the lemma below.

We choose the set $\mathcal{P}_n$ of [8] equal to $\{p_w : w \in B_n\}$, where $B_n \subset \ell^\infty(\mathcal{X})$ is the measurable set as in Theorem 2.1, with $C$ a large constant. In view of the first inequality of Lemma 3.1, for sufficiently large $n$, the $4\varepsilon_n$-entropy of $\mathcal{P}_n$ relative to the Hellinger distance is bounded above by the $3\varepsilon_n$-entropy of the set $B_n$ relative to the uniform distance, which is bounded by $6Cn\varepsilon_n^2$, by Theorem 2.1. This verifies (2.2) of [8]. The prior probability $\Pi(\mathcal{P}_n^c)$ outside the set $\mathcal{P}_n$, as in (2.3) of [8], is bounded by the probability of the event $\{W \notin B_n\}$, which is bounded by $e^{-Cn\varepsilon_n^2}$, by Theorem 2.1. Finally, by the second and third inequalities of Lemma 3.1, the prior probability as in (2.4) of [8], but with $\varepsilon_n$ replaced by a multiple of $\varepsilon_n$, is bounded below by the probability of the event $\{\|W - w_0\|_\infty < 2\varepsilon_n\}$, which is bounded below by $e^{-n\varepsilon_n^2}$, by Theorem 2.1. □

LEMMA 3.1. *For any measurable functions $v, w : \mathcal{X} \to \mathbb{R}$, we have the following:*

- $h(p_v, p_w) \leq \|v - w\|_\infty e^{\|v-w\|_\infty/2}$;
- $K(p_v, p_w) \lesssim \|v - w\|_\infty^2 e^{\|v-w\|_\infty}(1 + \|v - w\|_\infty)$;
- $V(p_v, p_w) \lesssim \|v - w\|_\infty^2 e^{\|v-w\|_\infty}(1 + \|v - w\|_\infty)^2$.

PROOF.  The triangle inequality and simple algebra give

$$h(p_v, p_w) = \left\| \frac{e^{v/2}}{\|e^{v/2}\|_2} - \frac{e^{w/2}}{\|e^{w/2}\|_2} \right\|_2 \leq 2 \frac{\|e^{v/2} - e^{w/2}\|_2}{\|e^{w/2}\|_2}.$$

Because $|e^{v/2} - e^{w/2}| = e^{w/2}|e^{v/2-w/2} - 1| \leq e^{w/2}e^{|v-w|/2}|v - w|/2$ for any $v, w \in \mathbb{R}$, the square of the right-hand side is bounded by

$$\frac{\int e^w e^{|v-w|}|v-w|^2 \, d\nu}{\int e^w \, d\nu} \leq e^{\|v-w\|_\infty}\|v-w\|_\infty^2.$$

This proves the first assertion of the lemma. We derive the other assertions from the first using the equivalence of $K$, $V$ and the Hellinger distance if the quotient of the densities is uniformly bounded. Because $w - \|v - w\|_\infty \leq v \leq w + \|v - w\|_\infty$, we have

$$\int e^w \, d\nu \, e^{-\|v-w\|_\infty} \leq \int e^v \, d\nu \leq \int e^w \, d\nu \, e^{\|v-w\|_\infty}.$$

Taking logarithms, we see that $-\|v - w\|_\infty \leq \log(\int e^v \, d\nu / \int e^w \, d\nu) \leq \|v - w\|_\infty$. Therefore,

$$\left\| \log \frac{p_v}{p_w} \right\|_\infty = \left\| v - w - \log \frac{\int e^v \, d\nu}{\int e^w \, d\nu} \right\|_\infty \leq 2\|v - w\|_\infty.$$

The second and third inequalities now follow from the first by Lemma 8 of [10].  □

3.2. *Classification.*  Suppose that we observe a random sample of vectors $(X_1, Y_1), \ldots, (X_n, Y_n)$ from the distribution of $(X, Y)$, where $Y$ takes its values in the set $\{0, 1\}$ and $X$ takes its values in some measurable space $(\mathcal{X}, \mathcal{A})$. Consider estimating the binary regression function $f_0(x) = \Pr(Y = 1 | X = x)$. Given a fixed, measurable function $\Psi \colon \mathcal{X} \to (0, 1)$, we may construct a prior on the set of regression functions as $f_W$ for a Gaussian process $W = (W_x : x \in \mathcal{X})$ and $f_w$ the function

$$f_w(x) = \Psi(w_x).$$

Here, $w_x$ denotes the value at $x$ of a function $w \colon \mathcal{X} \to \mathbb{R}$. The likelihood for $(X, Y)$ factorizes as

$$p_w(x, y) = f_w(x)^y(1 - f_w(x))^{1-y}g(x),$$

that is, into the conditional likelihood of $Y$ given $X$ and the marginal likelihood $g$ for $X$. As this causes the marginal density $g$ to cancel from the posterior distribution for $f_w$, it is not necessary to put a prior on $g$. We can set the distribution of $X$ equal to the "true" distribution $G$ into all of the following and incorporate it into the dominating measure $\nu$ so that the factor $g$ can be omitted from the likelihood. We assume that $f_0$ is never zero and then, with some abuse of notation, can define $w_0$ by the equation $f_0 = \Psi(w_0)$.

The link function $\Psi$ is assumed to be a known differentiable function with bounded derivative $\psi$. Link functions that lead to agreement between the metrics on the set of densities $p_w$ and the $L_2$-metric on the set of functions $w$ are especially attractive in the present context. The logistic link function qualifies in this respect. More generally, there is perfect agreement whenever the function $\psi/(\Psi(1-\Psi))$ is uniformly bounded.

Implementations of this model are described in [19, 20, 27, 34]. The first three follow [1] in defining latent variables and setting up an MCMC scheme. For probit regression, the latent variable is a single Gaussian variable $Z_i$ which models $Y_i = 1_{Z_i > 0}$ and, given $X_i = x$, possesses mean $w_x$ and variance 1. Logistic or other link functions are approximated by scale mixtures of Gaussian links, with an additional latent scale variable. [27] proposes to compute an approximation to the posterior distribution, either of Laplace form or by an algorithm termed "expectation propagation," both applicable to general priors.

THEOREM 3.2. (i) *Suppose that the function $\psi/(\Psi(1-\Psi))$ is bounded. Let $W$ be a Borel measurable, zero-mean, tight Gaussian random element in $L_2(G)$. Suppose that $w_0 = \Psi^{-1}(f_0)$ is contained in the support of $W$ and let $\phi_{w_0}$ be the function in (1.2) with $\|\cdot\|$ the $L_2(G)$-norm. Then, the posterior distribution relative to the prior $p_W$ satisfies $\mathrm{E}_0\Pi_n(w : \|p_w - p_0\|_{G,2} > M\varepsilon_n | X_1, Y_1, \ldots, X_n, Y_n) \to 0$ for any sufficiently large constant $M$ and $\varepsilon_n$ given by (1.3).*

(ii) *Suppose that the function $w_0 = \Psi^{-1}(f_0)$ is bounded. Let $W$ be a Borel measurable, zero-mean, tight Gaussian random element in $\ell^\infty(\mathcal{X})$. Suppose that $w_0$ is contained in the support of $W$ and let $\phi_{w_0}$ be the function in (1.2) with $\|\cdot\|$ the uniform norm. The same conclusion is then true.*

PROOF. This follows from combining Theorem 2.1 of [8] and Theorem 2.1 of the present paper, in the same fashion as Theorem 3.1 was proved by combining these two results. The details are as follows.

(i) Because the densities $p_w$ are uniformly bounded, the $L_2$-norm on the set of densities is bounded above by a multiple of the Hellinger distance. Thus, we can apply Theorem 2.1 of [8] with $d$ equal to the $L_2(G)$-norm. The square $L_2(G)$-norm and the Kullback–Leibler quantities $K$ and $V$ on the densities $p_w$ are all bounded above by multiples of the square $L_2(G)$-norm on the functions $w$, by Lemma 3.2 below. Therefore, Theorem 2.1, with $\|\cdot\|$ the $L_2(G)$-norm, allows us to bound the quantities in Theorem 2.1 of [8].

(ii) If the function $w_0 = \Psi^{-1}(f_0)$ is bounded, then so are the functions in a uniform neighborhood of it and so is the function $\psi/(\Psi(1 - \Psi))$ on the relevant domain. The proof can next be completed as before. □

The theorem can be extended to link functions with an unbounded function $\psi/(\Psi(1 - \Psi))$, even if the function $w_0 = \Psi^{-1}(f_0)$ is unbounded, by using appropriate norms on the Gaussian process. For instance, the probit link function can be treated as soon as the function $w_0 = \Psi^{-1}(f_0)$ is contained in $L_4(G)$, with a combination of the $L_2((w_0^2 \vee 1) \cdot G)$ and $L_4(G)$-norms on the Gaussian process. This can be proven in the same way as the preceding theorem, using Lemma 3.2 below.

For general link functions, the relationship between the appropriate norms on the densities and the norm on the Gaussian process is moderated by the function $S : \mathbb{R}^2 \to \mathbb{R}$ given by

$$S(w, w_0) = \sup_{v : v \in [w, w_0] \cup [w_0, w]} \left| \frac{\psi}{\Psi(1 - \Psi)}(v) \right| \vee 1.$$

LEMMA 3.2.  *If $\Psi$ possesses a bounded derivative $\psi$, then, for any measurable functions $v, w : \mathcal{X} \to \mathbb{R}$ and any $r > 1$, we have the following:*

- $\|p_v - p_w\|_r = 2^{1/r} \|\Psi(v) - \Psi(w)\|_{r,G} \le 2^{1/r} \|\psi\|_\infty \|v - w\|_{r,G}$;
- $K(p_w, p_{w_0}) \le \|(w - w_0)\sqrt{S(w, w_0)}\|_{2,G}^2$;
- $V(p_w, p_{w_0}) \le \|(w - w_0)S(w, w_0)\|_{2,G}^2$.

*For $\Psi$ the distribution function of the logistic distribution, the function $S$ is uniformly bounded. For $\Psi$ the distribution of the normal distribution, $K(p_w, p_{w_0})$ and $V(p_w, p_{w_0})$ are bounded above by a multiple of $\|w - w_0\|_{2,G_0}^2 + \|w - w_0\|_{G,4}^4$ where $G_0$ is the measure defined by $dG_0 = (w_0^2 \vee 1)\,dG$.*

PROOF.  The first assertion follows immediately from the fact that $|p_v(x, 0) - p_w(x, 0)| = |p_v(x, 1) - p_w(x, 1)| = |\Psi(v_x) - \Psi(w_x)|$ for any $x$. For the second inequality, we consider, for fixed $w_0 \in \mathbb{R}$, the function $g_{w_0} : \mathbb{R} \to \mathbb{R}$ given by

$$g_{w_0}(w) = \Psi(w_0) \log \frac{\Psi(w_0)}{\Psi(w)} + (1 - \Psi(w_0)) \log \frac{1 - \Psi(w_0)}{1 - \Psi(w)}.$$

The derivative of this function is $g'_{w_0}(w) = (\psi/\Psi(1 - \Psi))(w)(\Psi(w) - \Psi(w_0))$. In view of the definition of $S$ and Taylor's theorem, it follows that $|g'_{w_0}(w)| \le S(w, w_0)(w - w_0)^2$. The second assertion is then clear from the fact that $K(p_w, p_{w_0}) = \int g_{w_0}(w)\,dG$.

For the third inequality, we note that, by Taylor's theorem,

$$\left| \log \frac{\Psi(w)}{\Psi(w_0)} \right| \vee \left| \log \frac{1 - \Psi(w)}{1 - \Psi(w_0)} \right| \le S(w, w_0)|w - w_0|.$$

Since $V(p_w, p_{w_0})$ is a weighted integral of the squares of the quantities on the left-hand side, the third inequality follows.

For $\Psi$ the logistic distribution function, the function $\psi/(\Psi(1-\Psi))$ is easily seen to be bounded. The normal distribution function $\Psi$ satisfies $\psi/(\Psi(1-\Psi))(x) \asymp x$ as $x \to \pm\infty$ and is hence bounded by a multiple of $|x| \vee 1$. It follows that $|S(w, w_0)| \leq (|w_0| + |w_0 - w|) \vee 1$. Substituting this into the bounds on $K(p_w, p_{w_0})$ and $V(p_w, p_{w_0})$ readily yields the last assertion of the lemma. $\quad\square$

3.3. *Regression with fixed covariates.* Suppose that we observe independent variables $Y_1, \ldots, Y_n$ following the regression model $Y_i = w_0(x_i) + e_i$ for unobservable $N(0, \sigma_0^2)$-distributed errors $e_i$ and fixed, known elements $x_1, \ldots, x_n$ of a set $\mathcal{X}$. Consider estimating the regression function $w_0$.

As a prior on $w$, we use a Gaussian process $(W_x : x \in \mathcal{X})$. As this is conjugate, implementation is straightforward (see, e.g., [27], Chapter 2). If the standard deviation $\sigma$ of $e$ is not known, then we may also put a prior on $\sigma$, which we assume to be supported on a given interval $[a, b] \subset (0, \infty)$ with a Lebesgue density that is bounded away from zero. Unfortunately, the popular inverse Gamma prior does not satisfy the latter condition.

The natural semimetric for this problem is the $L_2(\mathbb{P}_n^x)$-norm for the empirical measure $\mathbb{P}_n^x = n^{-1} \sum_{i=1}^n \delta_{x_i}$ of the design variables. For fixed $n$, the Gaussian stochastic process $(W_x : x \in \mathcal{X})$ is important at the design points only and must be viewed as a Borel measurable map in the Banach space $L_2(\mathbb{P}_n^x)$. As this varies with $n$, it is more convenient to view it as a map in the space $\ell^\infty(\mathcal{X})$ of bounded functions on $\mathcal{X}$, whose norm is stronger than any of the $L_2(\mathbb{P}_n^x)$-norms.

THEOREM 3.3. *Let $W$ be a zero-mean, tight Gaussian random element in $\ell^\infty(\mathcal{X})$ and suppose that $w_0$ is contained in the support of $W$. Furthermore, let $S$ be a random variable with values in an interval $[a, b] \subset (0, \infty)$ that contains $\sigma_0$. Let $\phi_{w_0}$ be the function in (1.2) with $\|\cdot\|$ the supremum norm on $\ell^\infty(\mathcal{X})$. Then, the posterior distribution satisfies $\mathrm{E}_0 \Pi_n((w, \sigma) : \|w - w_0\|_n + |\sigma - \sigma_0| > M\varepsilon_n | Y_1, \ldots, Y_n) \to 0$ for any sufficiently large constant $M$ and $\varepsilon_n$ given by (1.3).*

PROOF. Let $\|\cdot\|_n$ be the $L_2(\mathbb{P}_n^x)$-norm. For the case where the prior on $\sigma$ is degenerate at the true value, it is shown in [11] that the rate of posterior contraction is faster than $\varepsilon_n$, for which there exist sets $\mathcal{W}_n$ satisfying

$$\log N(\varepsilon_n, \mathcal{W}_n, \|\cdot\|_n) \leq n\varepsilon_n^2,$$

$$\Pi_n(\mathcal{W}_n^c) \leq e^{-2n\varepsilon_n^2},$$

$$\Pi_n(w : \|w - w_0\|_n \leq \varepsilon_n) \geq e^{-n\varepsilon_n^2}.$$

This result is based on comparisons of the Kullback–Leibler divergence and variance to the square of the norm $\|\cdot\|_n$, and the construction of tests. It can be extended to the present case of an unknown scale that is bounded away from zero and infinity. The theorem then follows from Theorem 2.1. □

3.4. *White noise model.* Suppose that we observe a sample path of the stochastic process $X^{(n)} = (X_t^{(n)} : 0 \leq t \leq 1)$, defined structurally as, for a given function $w_0 \in L_2[0,1]$,

$$X_t^{(n)} = \int_0^t w_0(s)\,ds + \frac{1}{\sqrt{n}}B_t,$$

for a standard Brownian motion $B$. Consider estimating the function $w_0$. More formally, the statistical experiment consists of the set of induced distributions of the process $X^{(n)}$ on the Borel $\sigma$-field of the space $C[0,1]$ of continuous functions equipped with the uniform norm, as the parameter $w$ varies over a given subset of $L_2[0,1]$.

Consider the prior on the model obtained by modeling the parameter $w$ as the sample path of a Gaussian process $W$ with values in the space $L_2[0,1]$. As this is conjugate, the practical implementation is straightforward.

It is immediate from combining the preceding proposition with Theorem 2.1 that the rate of posterior contraction is determined by equations (1.2)–(1.3).

THEOREM 3.4. *Let $W$ be a zero-mean, tight Gaussian random element in $L_2[0,1]$ and suppose that $w_0$ is contained in the support of $W$. Let $\phi_{w_0}$ be the function in (1.2) with $\|\cdot\|$ the $L_2[0,1]$-norm. Then, the posterior distribution satisfies $\mathrm{E}_0\Pi_n(w : \|w - w_0\|_2 > M\varepsilon_n | X^{(n)}) \to 0$ for any sufficiently large constant $M$ and $\varepsilon_n$ given by (1.3).*

**4. Examples of Gaussian priors.** In this section, we give a number of examples of Gaussian process priors and compute their concentration functions (1.2) for "true parameters" of interest. We are especially interested in exhibiting processes that give the "correct" rates for true parameters of varying smoothness.

4.1. *Brownian motion and its primitives.* For modeling functions on the one-dimensional unit interval, Brownian motion is a good starting point. It can be viewed as a map into the space $C[0,1]$, but also as a map in $L_r[0,1]$. This does not affect its RKHS and small ball probabilities, which are both well known. The RKHS of Brownian motion is the collection of absolutely continuous functions $w : [0,1] \to \mathbb{R}$ with $w(0) = 0$ and $\int w'(t)^2\,dt < \infty$ with RKHS-norm $\|w\|_{\mathbb{H}} = \|w'\|_2$. The small ball probabilities of Brownian motion satisfy (cf. [24]), as $\varepsilon \downarrow 0$, for any $r \in [1, \infty]$,

$$-\log \Pr(\|W\|_r < \varepsilon) \asymp \left(\frac{1}{\varepsilon}\right)^2.$$

The support of Brownian motion in $C[0,1]$ is the set of all functions with $w(0) = 0$. Interestingly, the support as a map in $L_r[0,1]$, for $r < \infty$, is the full space $L_r[0,1]$.

The sample paths of Brownian motion are tied down to 0 at 0 and this, of course, remains the case for the functions in its RKHS. This can be relaxed by starting the process at an independent standard normal variable. The RKHS of "Brownian motion started at random" is the space of functions $w : [0,1] \to \mathbb{R}$ with $\int w'(t)^2 \, dt < \infty$ and with square RKHS-norm $\|w\|_{\mathbb{H}}^2 = w(0)^2 + \|w'\|_2^2$.

The small probability leads, by way of (1.3), to the restriction $\varepsilon_n^{-2} \le n\varepsilon_n^2$, equivalently $\varepsilon_n \ge n^{-1/4}$, on the rate of contraction.

The concentration function (1.2) further depends on the position of the true parameter $w_0$ relative to the RKHS. We may compute this contribution by approximation of $w_0$ through a kernel smoother. For $\phi_\sigma(x) = \sigma^{-1}\phi(x/\sigma)$ a smooth kernel, the convolution $w_0 * \phi_\sigma$ is contained in the RKHS and has uniform distance of the order $\sigma^\beta$ (as $\sigma \to 0$) to a function $w_0 \in C^\beta[0,1]$ that is Lipschitz of the order $\beta \in (0,1]$ and has square RKHS-norm $w_0 * \phi_\sigma(0)^2 + \|(w_0 * \phi_\sigma)'\|_2^2$ of the order $\sigma^{-(2-2\beta)}$ (see below). The choice $\sigma \asymp \varepsilon^{1/\beta}$ readily yields that

$$\inf\{\|h'\|_2^2 : \|h - w_0\|_\infty < \varepsilon\} \lesssim \varepsilon^{-(2-2\beta)/\beta}.$$

The concentration function (1.2) is the sum of this and the small ball exponent $(1/\varepsilon)^2$. For $\beta \ge 1/2$, the contribution of the small ball probability to the concentration function (1.2) dominates and the rate of contraction $\varepsilon_n$ is $n^{-1/4}$. For $\beta \in (0, 1/2)$, the contribution as in the preceding display dominates and will yield a rate not faster than $n^{-\beta/2}$. In particular, higher smoothness of the true parameter $w_0$ does not lead to a higher rate of contraction than $n^{-1/4}$ for $\beta > 1/2$.

For this reason, or by intuition, Brownian motion may be considered to be too rough as a prior. Integrating the sample paths one or more times will remedy this and will give Gaussian priors of smoothness $3/2, 5/2, \ldots$. To fill the gaps between these numbers, we consider, more generally, fractional integrals and fractional Brownian motion in the next sections. For ordinary integrals, the result is simpler and as follows.

Define $I_{0+}^1 f$ as the function $t \mapsto \int_0^t f(s) \, ds$ and $I_{0+}^k f$ as $I_{0+}^1(I_{0+}^{k-1} f)$.

**THEOREM 4.1.** *Let $W$ be a standard Brownian motion and $Z_0, \ldots, Z_k$ independent standard normal random variables. The RKHS of the process* $t \mapsto I_{0+}^k W_t + \sum_{i=0}^k Z_i t^i / i!$ *is the Sobolev space* $H^{k+1}[0,1]$ *with square norm* $\|h\|_{\mathbb{H}}^2 = \|h^{(k+1)}\|_2^2 + \sum_{i=0}^k h^{(i)}(0)^2$. *The concentration function of this process viewed as a map in* $C[0,1]$ *for an element* $w \in C^\beta[0,1]$ *for* $\beta \le k+1/2$ *satisfies* $\phi_w(\varepsilon) = O(\varepsilon^{-(2k-2\beta+2)/\beta})$ *as* $\varepsilon \downarrow 0$.

PROOF.    That the RKHS takes the present form is well known. See, for instance, [30] for a self-contained proof.

The concentration function of $I_{0+}^k W$ around zero satisfies $\phi_0(\varepsilon) \asymp \varepsilon^{-1/(k+1/2)}$ as $\varepsilon \to 0$ by Theorem 2.1 of [23]. The concentration function of the process $t \mapsto \sum_{i=0}^k Z_i t^k$ is of the order $\log(1/\varepsilon)$ and hence is much smaller. Thus, the concentration function of $I_{0+}^k W + \sum_{i=0}^k Z_i t^i / i!$ around zero is of the order $\varepsilon^{-1/(k+1/2)}$.

To compute the concentration function around a function $w \in C^\beta[0,1]$, we utilize convolutions $w * \phi_\sigma$ with a smooth higher-order kernel $\phi_\sigma$ with scale $\sigma$. As in kernel density estimation, the uniform distance between $w * \phi_\sigma$ and $w$ is of the order $\sigma^\beta$. The functions $w * \phi_\sigma$ belong to the RKHS. By writing $(w * \phi_\sigma)^{(l)} = w^{(\underline{\beta})} * \phi_\sigma^{(l-\underline{\beta})}$ for $\underline{\beta}$ the largest integer smaller than $\beta$, we see that $\|(w * \phi_\sigma)^{(l)}\|_\infty$ is bounded above by $\sigma^{-(l-\beta)}$ if $w \in C^\beta[0,1]$ and $l \ge \beta$ and hence that the RKHS-norm of $w * \phi_\sigma$ is of the order $\sigma^{-(2k-2\beta+2)}$ if $w \in C^\beta[0,1]$. Setting $\sigma \asymp \varepsilon^{1/\beta}$, we see that

$$\|w * \phi_\sigma - w\|_\infty \lesssim \varepsilon, \qquad \|w * \phi_\sigma\|_{\mathbb{H}}^2 \lesssim \varepsilon^{-(2k-2\beta+2)/\beta}.$$

Thus, the approximation part of $\phi_w(\varepsilon)$ is of the order $\varepsilon^{-(2k-2\beta+2)/\beta}$. For $\beta \le k+1/2$, this dominates the part $\varepsilon^{-1/(k+1/2)}$ resulting from the centered small ball probability.    □

For $\beta = k+1/2$, the concentration function, in the preceding theorem becomes $\phi_w(\varepsilon) = \varepsilon^{-1/\beta}$. For this function, inequality (1.3) is solved by

$$\varepsilon_n \asymp n^{-\beta/(2\beta+1)}.$$

This is the minimax rate for estimating a function that is known to be $\beta$-regular in various nonparametric models. Combination of the preceding theorem with the results on posterior contraction shows that the Gaussian prior in this case yields the optimal rate of convergence. For $\beta \ne k+1/2$, the Gaussian prior gives consistency with a rate, but the minimax rate is not achieved. This corresponds to an under- or over-smoothed prior.

Kimeldorf and Wahba [15] and Wahba [33] have considered priors of the type $t \mapsto \sqrt{b} I_{0+}^k W_t + \sqrt{a} \sum_{i=0}^k Z_i t^i / i$ in the setting of the regression model $Y_i = w(x_i) + e_i$. These priors are the same as in the preceding theorem, but with additional scaling factors $\sqrt{b}$ and $\sqrt{a}$. They show that if $a \to \infty$ and $b$ and $n$ are fixed, then the posterior mean of the regression function tends to the minimizer $w_n$ of the penalized least squares criterion

$$w \mapsto \frac{1}{n} \sum_{i=1}^{n} (w(x_i) - y_i)^2 + \frac{\sigma^2}{nb} \int_0^1 w^{(k)}(t)^2 \, dt,$$

where $\sigma^2$ is the variance of the regression errors $e_i$. Letting $a$ tend to infinity has the purpose of making the prior on the finite-dimensional, polynomial part diffuse, while the infinite-dimensional part of the prior is fixed. The preceding theorem considers the rate of contraction of the full prior as $n \to \infty$, for fixed $a$ and $b$, and hence is not directly comparable to the results of Kimeldorf and Wahba. However, some intriguing observations can be made. The penalized least squares estimator is well known to be a smoothing spline and is known to achieve the minimax rate $n^{-k/(2k+1)}$ for regression functions in $H^k[0,1]$ when the "smoothing parameter" $\lambda_n = \sigma^2/(nb)$ is set to satisfy $\lambda_n \asymp n^{-2k/(2k+1)}$. This would yield a scaling factor $b \asymp n^{-1/(2k+1)}$, meaning that the infinite-dimensional part of the prior would tend to zero. In contrast, the preceding theorem shows that a fixed value of $b$ yields a consistent posterior and a posterior achieving the optimal rate of contraction if the smoothness $\beta$ of the true parameter is equal to the smoothness $k + 1/2$ of the prior. (The theorem does not allow a diffuse prior on the polynomial part, but it can be checked that the theorem remains true if the Gaussians in this polynomial part have variance tending slowly to infinity.) It may be noted that the preceding theorem appears to indicate that the prior works best for functions in $H^{k+1/2}[0,1]$, not $H^k[0,1]$.

Wood and Kohn [34] implements the once-integrated Brownian motion prior within the setting of the binary regression model, where a large variance is used on the polynomial part.

4.2. *Riemann–Liouville process.* For $\alpha > 0$ and $W$ a standard Brownian motion, the *Riemann–Liouville process with Hurst parameter* $\alpha > 0$ is defined as

$$R_t = \int_0^t (t - s)^{\alpha - 1/2} \, dW_s, \qquad t \geq 0.$$

The process $R$ is a centered Gaussian process with continuous sample paths. It can be viewed as a multiple of the $(\alpha + 1/2)$-fractional integral of the "derivative $dW$ of Brownian motion." For $\alpha > 0$ and a (deterministic) measurable function $f$ on $[0,1]$, the (left-sided) *Riemann–Liouville fractional integral of $f$ of order $\alpha$* (if it exists) is defined as

$$I_{0+}^\alpha f(t) = \frac{1}{\Gamma(\alpha)} \int_0^t (t - s)^{\alpha - 1} f(s) \, ds.$$

For $\alpha$ a natural number, the function $I_{0+}^{\alpha} f$ is just the $\alpha$-fold iterated integral of $f$ and for $\alpha > 1/2$, the Riemann–Liouville process is equal to $\Gamma(\alpha + 1/2)I_{0+}^{\alpha - 1/2}W$ for $I_{0+}^{\alpha}$ the fractional integral. It can be shown that $I_{0+}^{\alpha}$ maps $\beta$-regular functions into $\alpha + \beta$-regular functions (if $\alpha + \beta$ is not an integer; see [28]). Since Brownian motion is "regular of order $1/2$," the Riemann–Liouville process $R$ is a good model for "$\alpha$-regular functions." This intuition is corroborated by the rate results in this section. For a proof of the following theorem, see Examples 9 and 15 in [13].

THEOREM 4.2. *The RKHS of the Riemann–Liouville process with parameter* $\alpha > 0$ *viewed as a random element in* $C[0, 1]$ *is* $\mathbb{H} = I_{0+}^{\alpha + 1/2}(L_2[0, 1])$ *and the RKHS-norm is given by*

$$\|I_{0+}^{\alpha + 1/2} f\|_{\mathbb{H}} = \frac{\|f\|_2}{\Gamma(\alpha + 1/2)}.$$

The Riemann–Liouville process is appropriate for approximating $C^{\alpha}$-functions, except that its definition as an integral from 0 means that its sample paths and their derivatives are tied down at zero. For $\alpha > 0$ and $\underline{\alpha}$ the biggest integer smaller than $\alpha$, we shall instead consider the process

$$(4.1) \qquad\qquad X_t = \sum_{k=0}^{\underline{\alpha}+1} Z_k t^k + R_t^{\alpha},$$

where $Z_1, \ldots, Z_{\underline{\alpha}+1}, R^{\alpha}$ are independent, $Z_i$ is standard normal and $R^{\alpha}$ is a Riemann–Liouville process with Hurst index $\alpha$. As before, we view this process as a random element in $C[0, 1]$.

THEOREM 4.3. *The support of the process* $X$ *is the whole space* $C[0, 1]$. *For any* $w \in C^{\alpha}[0, 1]$, *the concentration function of* $X$ *satisfies* $\phi_w(\varepsilon) = O(\varepsilon^{-1/\alpha})$ *as* $\varepsilon \downarrow 0$.

The proof of this theorem is deferred to Section 5. For $\alpha$ not an integer, it can be seen by inspection of the proof that the theorem remains true if $X$ is replaced by the process $\sum_{k=0}^{\underline{\alpha}} Z_k t^k + R_t^{\alpha}$.

For the concentration function $\phi_w(\varepsilon) = \varepsilon^{-1/\alpha}$, inequality (1.3) is solved by $\varepsilon_n = n^{-\alpha/(2\alpha+1)}$. This is the minimax rate for estimating a function that is known to be $\alpha$-regular in various nonparametric models. Combination of the preceding theorem with the results on posterior contraction therefore shows that the Gaussian prior (4.1) yields the optimal rate of convergence in various settings. This is true, for instance, in the settings of density estimation using a prior of the form $t \mapsto ce^{X_t}$ on the density, Gaussian regression using $X$ as a prior regression function and classification using a prior $\Psi(X_t)$ on the probability $\Pr(Y = 1 | X = t)$.

4.3. *Fractional Brownian motion.* Fractional Brownian motion offers a different starting point for constructing a Gaussian process of a given smoothness level. By definition, *fractional Brownian motion (fBm)* with Hurst parameter $\alpha \in (0, 1)$ is the zero-mean Gaussian process $X = (X_t : t \in [0, 1])$ with continuous sample paths and covariance function

$$\mathrm{E} X_s X_t = \tfrac{1}{2}(s^{2\alpha} + t^{2\alpha} - |t - s|^{2\alpha}).$$

The choice $\alpha = 1/2$ yields ordinary Brownian motion. To obtain a process of a given smoothness $\alpha > 1$, we can take ordinary integrals of fractional Brownian motion.

The conclusions using fractional Brownian motion are the same as for the Riemann–Liouville process.

THEOREM 4.4. *Consider the fractional Brownian motion with Hurst parameter $\alpha \in (0, 1)$ as a random element in $C[0, 1]$. For $w \in C^\alpha[0, 1]$ with $w(0) = 0$, we have $\phi_w(\varepsilon) = O(\varepsilon^{-1/\alpha})$ as $\varepsilon \to 0$.*

PROOF. For the fBm $X$, we have the representation

$$X_t = c_\alpha \int_{-\infty}^\infty ((t - s)_+^{\alpha - 1/2} - (-s)_+^{\alpha - 1/2}) \, dW_s,$$

where $W$ is a double-sided Wiener process and $c_\alpha$ a positive constant [25]. In other words, we have $X = c_\alpha R + c_\alpha Z$, where $R$ and $Z$ are independent processes, $R$ is a RL-process with parameter $\alpha$ and $Z$ is defined by

$$Z_t = \int_{-\infty}^0 ((t - s)^{\alpha - 1/2} - (-s)^{\alpha - 1/2}) \, dW_s.$$

By Lemma 3.2 of [23], $-\log \Pr(\|Z\| < \varepsilon) = o(\varepsilon^{-1/\alpha})$ as $\varepsilon \to 0$, with $\|\cdot\|$ the supremum norm, hence also with $\|\cdot\|$ the $L_2$-norm. The theorem therefore follows from the results for the RL-process and Theorem 2.3. □

4.4. *Truncated series.* Any Gaussian variable in a separable Banach space can be expanded as an infinite series $\sum_i Z_i h_i$ for i.i.d. standard normal variables $Z_i$ and elements $h_i$ from its RKHS. By Theorem 2.2, the prior obtained by truncating the series at a sufficiently high level will have the same concentration function and will hence lead to the same posterior rate of contraction. Because finite sums may be easier to handle, it is interesting to investigate special expansions and the numbers of terms that need to be retained in order to obtain the same contraction rate. In this section, we consider this question for fractional Brownian motion.

By Theorem 4.4, the fBm with Hurst parameter $\alpha \in (0, 1)$ as a prior for a true signal $w_0$ which is Hölder continuous of order $\alpha$ leads to a concentration function satisfying $\phi_{w_0}(\varepsilon) \lesssim \varepsilon^{-1/\alpha}$. Consequently, inequality (1.3) is satisfied for $\varepsilon_n \asymp n^{-\alpha/(1+2\alpha)}$. This implies, for instance, that the rate of posterior contraction in the white noise model (see Theorem 3.4) is equal to the minimax rate relative to the $L_2$-norm.

By Theorem 2.2, for any series expansion $X = \sum_k Z_k h_k$ of the fBm $X$, the truncated series $X^K = \sum_{k=1}^K Z_k h_k$ also gives the optimal rate of posterior contraction if

$$(4.2) \qquad 10\,\mathrm{E}\|X^K - X\|_2^2 \le \frac{1}{n}.$$

It is known that for any such expansion of the fBm, the truncated series $X^K$ satisfies (cf. [17])

$$\mathrm{E}\|X - X^K\|_2 \gtrsim \left(\frac{1}{K}\right)^\alpha.$$

A given expansion is therefore called *rate-optimal* (for the $L_2$-norm) if $\mathrm{E}\|X - X^K\|_2 \lesssim K^{-\alpha}$. Several explicit rate-optimal expansions of the fBm are known (see, e.g., [2, 6, 7, 14]). For these rate-optimal expansions, (4.2) is fulfilled as soon as the number of terms in the expansion satisfies $K = K_n \asymp Cn^{1/(2\alpha)}$ for a large constant $C$. This is somewhat larger than the dimension $n^{1/(2\alpha+1)}$ found in the following section, which also arises in the usual bias-variance trade-off of series estimators.

4.5. *Finite sums.* Replacing the coefficients in a series expansion by Gaussian variables is a natural method to construct a Gaussian prior on a set of functions. In this section, we use truncated series and study the effect of varying the variances of the Gaussian variables.

Series priors have been implemented in [20] in the density estimation model of Section 3.1 and in [21] in semiparametric regression, with Fourier-type series and coefficients with exponentially decreasing variances. The priors are easy to implement in Gaussian regression.

Because wavelet expansions give easy control of various norms, we consider here expansions

$$w = \sum_{j=1}^\infty \sum_{k=1}^{2^{jd}} w_{j,k} \psi_{j,k}$$

of functions $w: [0,1]^d \to \mathbb{R}$ on a double-indexed basis $\{\psi_{j,k} : j = 1, 2, \ldots, k = 1, \ldots, 2^{jd}\}$ of bounded functions $\psi_{j,k}: [0,1]^d \to \mathbb{R}$. (The unit cube could be replaced by another compact subset of $\mathbb{R}^d$.) We consider these functions with the norms

$$\|w\|_2 = \sum_{j=1}^\infty \left(\sum_{1 \le k \le 2^{jd}} |w_{j,k}|^2\right)^{1/2},$$

$$\|w\|_\infty = \sum_{j=1}^\infty 2^{jd/2} \max_{1\le k\le 2^{jd}} |w_{j,k}|,$$

$$\|w\|_{\beta|\infty,\infty} = \sup_{1\le j<\infty} 2^{j\beta} 2^{jd/2} \max_{1\le k\le 2^{jd}} |w_{j,k}|.$$

For the base functions $\psi_{j,k}$ derived from suitable orthonormal wavelets in $L_2[0,1]^d$, these norms correspond to the $L_2$-norm, the supremum norm and the Besov $(\beta,\infty,\infty)$-norm, respectively. The last norm measures smoothness of order $\beta$, weaker than a Hölder norm of the same order.

For given truncation levels $J_\alpha$, which will tend to infinity with $n$, we consider a Gaussian prior of the type

$$(4.3) \qquad W = \sum_{j=1}^{J_\alpha} \sum_{k=1}^{2^{jd}} \mu_j Z_{j,k} \psi_{j,k},$$

where the $\mu_j$ are positive numbers and the $Z_{j,k}$ are i.i.d. standard normal variables. The number of terms in the random series is $O(2^{J_\alpha d})$. For a transparent description of the main results, we set this number equal to the integer closest to the solution $\bar{J}_\alpha$ of the equation, for a given $\alpha > 0$,

$$2^{\bar{J}_\alpha d} = n^{d/(2\alpha+d)}.$$

This dimension is well known to be the optimal dimension of a finite-dimensional model if the true parameter is known to be regular of order $\alpha$. We next study the rate of posterior contraction if the true parameter is $\beta$-regular under a variety of choices of the coefficients $\mu_j$ and for a general $\beta > 0$, which may be smaller or larger than the "nominal" value $\alpha$.

The contribution $W_j = \sum_{k=1}^{2^{jd}} \mu_j Z_{j,k} \psi_{j,k}$ of the $j$th level to the prior satisfies

$$\mathrm{E}\|W_j\|_2^2 = \mu_j^2 2^{jd}.$$

Therefore, the choice $\mu_j = 2^{-jd/2}$ gives all levels the same amount of prior uncertainty. It is natural to choose the constants $\mu_j$ so that the numbers $2^{jd/2}\mu_j$ are nonincreasing, but we shall allow these numbers to tend to zero as $j \to \infty$. If $2^{jd/2}\mu_j \to 0$, then the higher levels receive less weight and hence the prior tends to be of lower dimension than the nominal dimension $2^{J_\alpha d}$. This may be advantageous if the true parameter is of higher regularity (i.e., $\beta > \alpha$), for which the optimal dimension $2^{J_\beta d}$ is indeed smaller. On the other hand, if the true parameter is less regular (i.e., $\beta < \alpha$), then the nominal dimension $2^{J_\alpha d}$ is already too small and this would be exacerbated by putting lower weight on the higher levels. We shall show that the choice $2^{jd/2}\mu_j = 2^{-j\beta}$ is a good compromise: it yields the optimal rate of contraction $n^{-\beta/(2\beta+d)}$ if $\beta \ge \alpha$ and the "optimal rate using a $2^{J_\alpha d}$-dimensional model" $n^{-\beta/(2\alpha+d)}$ if $\beta \leq \alpha$. The choice $2^{jd/2}\mu_j = 1$, which gives equal weight to all levels, is no worse than this if $\beta \leq \alpha$, but yields only the rate $n^{-\alpha/(2\alpha+d)}$ for $\beta \geq \alpha$.

To be precise, in the following, we establish these rates up to logarithmic factors. The proof of the following theorem can be found in Section 5.

THEOREM 4.5. *Let $W$ be the Gaussian process given by* (4.3) *viewed as a map in $\ell^\infty[0,1]^d$, with $\mu_j 2^{jd/2} = 2^{-ja}$ for some $a \geq 0$. Let $w_0 : [0,1]^d \to \mathbb{R}$ satisfy $\|w_0\|_{\beta|\infty,\infty} < \infty$. Then, for*

$$\varepsilon_n \geq \begin{cases} n^{-\beta/(2\alpha+d)}\log n, & \text{if } a \leq \beta \leq \alpha, \\ n^{-\alpha/(2\alpha+d)}\log n, & \text{if } a \leq \alpha \leq \beta, \\ n^{-a/(2a+d)}(\log n)^{d/(2a+d)}, & \text{if } \alpha \leq a \leq \beta, \\ n^{-\beta/(2a+d)}(\log n)^{d/(2a+d)}, & \text{if } \alpha \leq \beta \leq a, \end{cases}$$

*there exists a measurable set $B_n \subset \ell^\infty[0,1]^d$ such that*

$$(4.4) \qquad\qquad \log N(3\varepsilon_n, B_n, \|\cdot\|_\infty) \leq 6Cn\varepsilon_n^2,$$

$$(4.5) \qquad\qquad \Pr(W \notin B_n) \leq e^{-Cn\varepsilon_n^2},$$

$$(4.6) \qquad\qquad \Pr(\|W - w_0\|_\infty < 4\varepsilon_n) \geq e^{-n\varepsilon_n^2}.$$

## 5. Proofs.

PROOF OF THEOREM 2.1. Inequality (2.4) is an immediate consequence of (1.3) and (4.16) of [16]. We need to prove existence of the sets $B_n$ such that the first and second inequalities in the theorem hold.

For $\mathbb{B}_1$ and $\mathbb{H}_1$ the unit balls in the Banach space $\mathbb{B}$ and the RKHS $\mathbb{H}$, respectively, and $M_n$ a positive constant, set

$$B_n = \varepsilon_n \mathbb{B}_1 + M_n \mathbb{H}_1.$$

By Borell's inequality ([3], Theorem 3.1), it follows that

$$\Pr(W \notin B_n) \leq 1 - \Phi(\alpha_n + M_n)$$

for $\Phi$ the distribution function of the standard normal distribution and $\alpha_n$ determined by

$$\Phi(\alpha_n) = \Pr(W \in \varepsilon_n \mathbb{B}_1) = e^{-\phi_0(\varepsilon_n)}.$$

For $C > 1$, set

$$M_n = -2\Phi^{-1}(e^{-Cn\varepsilon_n^2}).$$

Because $\phi_0(\varepsilon_n) \leq \phi_{w_0}(\varepsilon_n) \leq n\varepsilon_n^2$ by assumption (1.3), and $C > 1$, we have that $\alpha_n \geq -\frac{1}{2}M_n$, whence $\alpha_n + M_n \geq \frac{1}{2}M_n$ and

$$\Pr(W \notin B_n) \leq 1 - \Phi(\tfrac{1}{2}M_n) = e^{-Cn\varepsilon_n^2}.$$

We conclude that inequality (2.3) is satisfied.

If $h_1, \ldots, h_N$ are contained in $M_n \mathbb{H}_1$ and are $2\varepsilon_n$-separated for the norm $\| \cdot \|$, then the $\| \cdot \|$-balls $h_j + \varepsilon_n \mathbb{B}_1$ of radius $\varepsilon_n$ around these points are disjoint and hence

$$1 \geq \sum_{j=1}^{N} \Pr(W \in h_j + \varepsilon_n \mathbb{B}_1)$$

$$\geq \sum_{j=1}^{N} e^{(-1/2)\|h_j\|_{\mathbb{H}}^2} \Pr(W \in \varepsilon_n \mathbb{B}_1)$$

$$\geq N e^{-(1/2)M_n^2} e^{-\phi_0(\varepsilon_n)},$$

where the second inequality follows from (4.16) of [16]. If the $2\varepsilon_n$-net $h_1, \ldots, h_N$ is maximal in the set $M_n \mathbb{H}_1$, then the balls $h_j + 2\varepsilon_n \mathbb{B}_1$ cover $M_n \mathbb{H}_1$. It follows that

$$N(2\varepsilon_n, M_n \mathbb{H}_1, \| \cdot \|) \leq N \leq e^{(1/2)M_n^2} e^{\phi_0(\varepsilon_n)}.$$

By its definition, any point of the set $B_n$ is within distance $\varepsilon_n$ of some point of $M_n \mathbb{H}_1$. This implies that

$$\log N(3\varepsilon_n, B_n, \| \cdot \|) \leq \log N(2\varepsilon_n, M_n \mathbb{H}_1, \| \cdot \|)$$

$$\leq \tfrac{1}{2} M_n^2 + \phi_0(\varepsilon_n)$$

$$\leq 5Cn\varepsilon_n^2 + \phi_0(\varepsilon_n),$$

by the definition of $M_n$ if $e^{-Cn\varepsilon_n^2} < 1/2$, because $\Phi^{-1}(y) \geq -\sqrt{5/2 \log(1/y)}$ and is negative for every $y \in (0, 1/2)$. Since $\phi_0(\varepsilon_n) \leq \phi_{w_0}(\varepsilon_n) \leq n\varepsilon_n^2$, this concludes the verification of (2.2). $\quad \square$

Proof of Theorem 2.2. As a consequence of Borell's inequality (cf. [32], Proposition A2.1), we have

$$\Pr(\|W^n - W\| \geq \varepsilon_n) \leq 2e^{-\varepsilon_n^2/8\mathrm{E}\|W^n - W\|^2} \leq 2e^{-n\varepsilon_n^2/(8/10)}$$

since $\mathrm{E}\|W^n - W\|^2 \leq 1/(10n)$, by assumption. For $\varepsilon_n$ satisfying $n\varepsilon_n^2 \geq 4\log 4$, the right-hand side is bounded above by $\tfrac{1}{2} e^{-n\varepsilon_n^2}$. Because $\Pr(\|W^n - w_0\| < 3\varepsilon_n) \geq \Pr(\|W - w_0\| < 2\varepsilon_n) - \Pr(\|W^n - W\| \geq \varepsilon_n)$, it follows from (1.3) that

$$\Pr(\|W^n - w_0\| < 3\varepsilon_n) \geq \tfrac{1}{2} e^{-n\varepsilon_n^2} \geq e^{-4n\varepsilon_n^2}.$$

This completes the verification of (2.4) with $\varepsilon_n$ replaced by $2\varepsilon_n$.

We choose $B_n = 2\varepsilon_n \mathbb{B}_1 + M_n \mathbb{H}_1^n$ for $\mathbb{H}_1^n$ the unit ball of the RKHS associated to $W^n$, and $M_n = -2\Phi^{-1}(e^{-Cn\varepsilon_n^2})$, as in the proof of Theorem 2.1. Similarly to the observation in the preceding paragraph, we have

$$e^{-\phi_0^n(2\varepsilon_n)} := \Pr(\|W^n\| < 2\varepsilon_n) \geq \tfrac{1}{2} e^{-n\varepsilon_n^2} \geq e^{-4n\varepsilon_n^2}.$$

The verification of (2.3)–(2.4) with $2\varepsilon_n$ instead of $\varepsilon_n$ can now proceed exactly as in the proof of Theorem 2.1. For the first, we use that $C > 4$, so that again $\alpha_n = \Phi^{-1}(e^{-\phi_0^n(2\varepsilon_n)}) \geq \Phi^{-1}(e^{-4n\varepsilon_n^2}) \geq -\frac{1}{2}M_n$. For the second, we substitute the inequality $\phi_0^n(2\varepsilon_n) \leq 4n\varepsilon_n^2$.  □

PROOF OF THEOREM 2.3.   If $\|W^i - w^i\| < \varepsilon$ for every $i$, then $\|W - w\| < \varepsilon|I|$, where $|I|$ is the cardinality of $I$. Combined with the independence of the processes $W^i$, this implies that $\Pr(\|W - w\| < \varepsilon|I|) \geq \prod_i \Pr(\|W^i - w^i\| < \varepsilon)$. In view of Theorem 2 of [16], the concentration function $\phi_w(\varepsilon)$ of $W$ is bounded above by twice the negative logarithm of the left-hand side, which is bounded above by $2\sum_i \phi_{w^i}^i(\varepsilon)$, again by Theorem 2 of [16].  □

PROOF OF THEOREM 2.4.   It is easy to see that the RKHS $\mathbb{H}^a$ of the process $aW$ for a fixed value of $a$ is equal to the RKHS $\mathbb{H}$ of $W$, but with norm $\|h\|_{\mathbb{H}^a} = a^{-1}\|h\|_{\mathbb{H}}$. We define $B_n = K\varepsilon_n\mathbb{B}_1 + KM_n\mathbb{H}_1 = KB_n^1$ for $B_n^1$ the set $B_n$ appearing in Theorem 2.1. Because $A \leq K$ and $B_n$ is a cone, it is clear that $\Pr(AW \notin B_n) \leq \Pr(W \notin B_n^1) \leq e^{-Cn\varepsilon_n^2}$, by Theorem 2.1. Also, $N(3K\varepsilon_n, B_n, \|\cdot\|) \leq N(3\varepsilon_n, B_n^1, \|\cdot\|) \leq 6Cn\varepsilon_n^2$, again by Theorem 2.1.

By Theorem 2 of [16], for any fixed $a$ and $\varepsilon > 0$,

$$-\log\Pr(\|aW - w_0\| < 2\varepsilon)$$
$$\leq \inf\{\|h\|_{\mathbb{H}^a}^2 : \|h - w\| < \varepsilon\} - \log\Pr(\|aW\| < \varepsilon)$$
$$\leq \frac{1}{a^2}\inf\{\|h\|_{\mathbb{H}}^2 : \|h - w\| < \varepsilon\} - \log\Pr(\|W\| < \varepsilon/K)$$
$$\leq \frac{1}{k^2}\phi_w(\varepsilon/K)$$

for $a > k$ and $0 < k < 1 < K$. We apply this with $\varepsilon/K = \varepsilon_n$ and then apply (1.3) to arrive at (2.7).  □

For the proof of Theorem 4.3, we first recall some facts from fractional calculus, which can be found in [28].

Using Fubini's theorem, it can be seen that the fractional integration operators have the semigroup property $I_{0+}^\alpha I_{0+}^\beta = I_{0+}^{\alpha+\beta}$. The fractional integration operator acts on power functions as one would expect: for $\alpha > 0$, $\beta > -1$ and $f(t) = t^\beta$,

$$I_{0+}^\alpha f(t) = \frac{\Gamma(\beta+1)}{\Gamma(\alpha+\beta+1)}t^{\alpha+\beta}.$$

For $\alpha \in (0,1)$, the (left-sided) *Riemann–Liouville fractional derivative of $f$ of order $\alpha$* is defined by

$$D_{0+}^\alpha f(t) = \frac{1}{\Gamma(1-\alpha)}\frac{d}{dt}\int_0^t (t-s)^{-\alpha}f(s)\,ds = \frac{d}{dt}I_{0+}^{1-\alpha}f(t),$$

provided it exists. To define the fractional derivative for $\alpha \geq 1$, we introduce the notation $[\alpha]$ and $\{\alpha\}$ for the integer and fractional parts of $\alpha$, respectively. For general $\alpha > 0$, we define

$$D_{0+}^{\alpha} f = \left(\frac{d}{dt}\right)^{[\alpha]} D_{0+}^{\{\alpha\}} f.$$

In particular, $D_{0+}^{\alpha} f$ is just the $\alpha$th derivative of $f$ if $\alpha$ is an integer. Observe that $D_{0+}^{\alpha} f$ equals the $n$th derivative of $I_{0+}^{n-\alpha} f$, provided it exists, with $n = \underline{\alpha} + 1$. We say that $f$ has a *summable fractional derivative* $D_{0+}^{\alpha} f$ if $I_{0+}^{n-\alpha} f$ has $n-1$ continuous derivatives and the $(n-1)$th derivative is only absolutely continuous rather than differentiable.

Fractional integration and differentiation are inverse operations, in the sense that $D_{0+}^{\alpha} I_{0+}^{\alpha} = \mathrm{Id}$. However, in general, $D_{0+}^{\alpha}$ is not the right inverse of $I_{0+}^{\alpha}$. If $f \in L_1$ has a summable derivative of order $\alpha > 0$, then, with $n = [\alpha] + 1$,

$$I_{0+}^{\alpha} D_{0+}^{\alpha} f(t) = f(t) - \sum_{k=0}^{n-1} \frac{D_{0+}^{n-k-1}(I_{0+}^{n-\alpha} f)(0)}{\Gamma(\alpha - k)} t^{\alpha-k-1}.$$

LEMMA 5.1. *Suppose that $f$ is twice continuously differentiable and $f(0) = 0$. For $\alpha \in (1, 2)$, the function $f$ has a summable fractional derivative $D_{0+}^{\alpha} f$ and can be written as $f = I_{0+}^{\alpha} D_{0+}^{\alpha} f$. Furthermore,*

$$D_{0+}^{\alpha} f(t) = \frac{f'(0)}{\Gamma(2-\alpha)} t^{1-\alpha} + I_{0+}^{2-\alpha} f''(t).$$

PROOF. Since $f(0) = 0$, we have $f(t) = f'(0)t + I_{0+}^2 f''(t)$, whence

$$I_{0+}^{2-\alpha} f(t) = \frac{f'(0)}{\Gamma(4-\alpha)} t^{3-\alpha} + I_{0+}^{4-\alpha} f''(t).$$

Differentiating this twice and using the identity $\Gamma(1+x) = x\Gamma(x)$ yields the expression for $D_{0+}^{\alpha} f$. The formula preceding the lemma gives

$$f(t) = I_{0+}^{\alpha} D_{0+}^{\alpha} f(t) + \sum_{k=0}^{1} \frac{D_{0+}^{1-k}(I_{0+}^{2-\alpha} f)(0)}{\Gamma(\alpha - k)} t^{\alpha-k-1}.$$

For $k = 1$, we get $I_{0+}^{2-\alpha} f(0)$ in the numerator, which vanishes since $f$ is continuous. For $k = 0$, we get $D_{0+}^{\alpha-1} f(0)$. But since $f(0) = 0$, we have $f = I_{0+}^1 f'$ and hence $D_{0+}^{\alpha-1} f(0) = I_{0+}^{2-\alpha} f'(0) = 0$. □

Roughly speaking, fractional integration of order $\alpha$ improves the smoothness of a function by $\alpha$. More precisely, for $\lambda \in [0, 1]$ and $\alpha \in (0, 1)$ such that $\lambda + \alpha \neq 1$, it holds that $I_{0+}^{\alpha} : C_0^{\lambda}[0, 1] \to C^{\alpha+\lambda}[0, 1]$, where $C_0^{\lambda}[0, 1]$ are the functions $f \in C^\lambda[0,1]$ with $f(0) = 0$. The analogous assertion is true for the fractional derivative. If $0 < \alpha < \lambda \le 1$, then $C_0^\lambda[0,1] \subset I_{0+}^\alpha(L_1[0,1])$. On the space $I_{0+}^\alpha(L_1[0,1])$, the Riemann–Liouville fractional derivative $D_{0+}^\alpha$ coincides with the so-called *Marchaud fractional derivative* of order $\alpha$. The latter maps $C_0^\lambda[0,1]$ into $C^{\lambda-\alpha}[0,1]$. Hence, for $0 < \alpha < \lambda \le 1$, we have $D_{0+}^\alpha : C_0^\lambda[0,1] \to C^{\lambda-\alpha}[0,1]$.

In the following lemma, we use the customary notation $I_{0+}^\alpha = D_{0+}^{-\alpha}$ for $\alpha < 0$.

LEMMA 5.2. *Let $\lambda \in [0,1]$ and $\alpha \in [0,1)$ be such that $\alpha + \lambda \in (0,2)$ and $\alpha + \lambda \ne 1$. If $f \in C^\lambda[0,1]$ and $g \in L_1(\mathbb{R})$ has compact support and satisfies $\int g(u)\,du = 0$ and, in the case that $\alpha + \lambda > 1$, also $\int u g(u)\,du = 0$, then*

$$\|I_{0+}^\alpha(f * g)\|_\infty \lesssim \int |u|^{\alpha+\lambda} |g(u)|\,du.$$

PROOF. The conditions on $g$ imply that

$$(f * g)(s) = \int (f(s-u) - f(s)) g(u)\,du \qquad \text{for } s \in (0,1)$$

and we may assume that $f(0) = 0$. A change of variables shows that for $u \in \mathbb{R}$ and $t \in (0,1)$, we have

$$\frac{1}{\Gamma(\delta)} \int_0^t (t-s)^{\alpha-1} f(s-u)\,ds = I_{0+}^\alpha f(t-u),$$

the right-hand side vanishing by definition if $u > t$. Using the fact that $g$ has compact support to justify the interchanging of integrals, it follows that

$$(5.1) \qquad I_{0+}^\alpha(f * g)(t) = \int ((I_{0+}^\alpha f)(t-u) - (I_{0+}^\alpha f)(t)) g(u)\,du.$$

Because $I_{0+}^\alpha : C_0^\lambda[0,1] \to C^{\alpha+\lambda}[0,1]$, we have, for $\alpha + \lambda > 1$,

$$|(I_{0+}^\delta f)(t-u) - (I_{0+}^\delta f)(t) + u(I_{0+}^\alpha f)'(t)| \lesssim |u|^{\alpha+\lambda}.$$

Inserting this in the preceding display completes the proof in this case. If $\alpha + \lambda < 1$, then the preceding display is satisfied with the factor $u(I_{0+}^\alpha f)'(t)$ omitted and the proof is completed as before. $\square$

PROOF OF THEOREM 4.3. Let $Z = X - R^\alpha$ be the polynomial part of $X$ given in (4.1).

By Theorem 2.1 of [23], $-\log \Pr(\|R^\alpha\|_\infty < \varepsilon)$ behaves as a constant times $\varepsilon^{-1/\alpha}$ as $\varepsilon \to 0$. Because each of the probabilities $\Pr(\|Z_k t^k\|_\infty < \varepsilon)$ behaves as a constant times $\varepsilon$ as $\varepsilon \to 0$, $-\log \Pr(\|Z\|_\infty < \varepsilon)$ is bounded above by a constant times $\log(1/\varepsilon)$, which is much smaller than $\varepsilon^{-1/\alpha}$.

In view of Theorem 2.3, the concentration function $\phi_w(2\varepsilon)$ of $X$ is bounded by a multiple of the sum $\phi_{w-P}(2\varepsilon; R^\alpha) + \phi_P(2\varepsilon; Z)$ of the concentration functions of $R^\alpha$ and $Z$, where $w = w - P + P$ may be an arbitrary split. The RKHS of the process $Z$ is the set of polynomials $P_\xi = \sum_{i=0}^{\alpha+1} \xi_i t^i$ with square norm $\|P_\xi\|_{\mathbb{H}}^2 = \sum_{i=1}^{\alpha+1} \xi_i^2$. Therefore, for any such polynomial,

$$\phi_{P_\xi}(\varepsilon; Z) \lesssim \sum_{i=1}^{\underline{\alpha}+1} \xi_i^2 + \log(1/\varepsilon).$$

We shall apply this with polynomials such that $\phi_{w-P}(2\varepsilon; R^\alpha)$ becomes suitably small.

Let $\phi$ be a smooth, compactly supported, order-$\underline{\alpha} \vee 2$ kernel and, for $\sigma > 0$, define $\phi_\sigma(t) = \sigma^{-1}\phi(t/\sigma)$. We note that, automatically, $\int \phi'(t)\,dt = \int \phi''(t)\,dt = \int t\phi''(t)\,dt = 0$. Since $w \in C^\alpha$, we have $\|w - w * \phi_\sigma\|_\infty \lesssim \sigma^\alpha$, whence $\|w - w * \phi_\sigma\|_\infty \le \varepsilon$ if $\sigma = C\varepsilon^{1/\alpha}$ for an appropriate constant $C$.

Let $\gamma = \{\alpha\} \in (0,1]$ be the fractional part of $\alpha$. We first consider the case $\gamma \in (0, 1/2]$. By Taylor's theorem,

$$w * \phi_\sigma(t) = \sum_{k=0}^{\underline{\alpha}} \frac{(w^{(k)} * \phi_\sigma)(0)}{k!} t^k + I_{0+}^{\underline{\alpha}+1}(w^{(\underline{\alpha})} * \phi_\sigma')$$

$$= \sum_{k=0}^{\underline{\alpha}} \frac{(w^{(k)} * \phi_\sigma)(0)}{k!} t^k + I_{0+}^{\alpha+1/2} I_{0+}^{1/2-\gamma}(w^{(\underline{\alpha})} * \phi_\sigma'),$$

by the semigroup property of fractional integrals. The first function on the right is a polynomial $P_\sigma$ and the sum of squares of its coefficients can be seen to be bounded for small $\sigma$. By Theorem 4.2, the second function on the right belongs to the RKHS of $R^\alpha$ and has squared RKHS-norm equal to $\|I_{0+}^{1/2-\gamma}(w^{(\underline{\alpha})} * \phi_\sigma')\|_2^2$, which is $O(\sigma^{-1}) = O(\varepsilon^{-1/\alpha})$, by Lemma 5.2. We now split $w = w - P_\sigma + P_\sigma$ and approximate $w - P_\sigma$ with the function $w * \phi_\sigma - P_\sigma$ in the RKHS of $R^\alpha$.

In the case where $\gamma \in (1/2, 1]$, we apply Lemma 5.1, with the $\alpha$ and function $f$ in the lemma taken equal to the present $\gamma + 1/2$ and $I_{0+}^1(w^{(\underline{\alpha})} * \phi_\sigma')$, to obtain that

$$I_{0+}^1(w^{(\underline{\alpha})} * \phi_\sigma') = I_{0+}^{\gamma+1/2}(g_1 + g_2),$$

with

$$g_1(t) = \frac{(w^{(\underline{\alpha})} * \phi_\sigma')(0)}{\Gamma(3/2 - \gamma)} t^{1/2-\gamma},$$

$$g_2(t) = I_{0+}^{3/2-\gamma}(w^{(\underline{\alpha})} * \phi_\sigma'')(t).$$

By integrating the penultimate display $\underline{\alpha}$ times, we obtain $I_{0+}^{\underline{\alpha}+1}(w^{\underline{\alpha}} * \phi_\sigma') = I_{0+}^{\alpha+1/2}(g_1 + g_2)$ and hence, by Taylor's theorem,

$$w * \phi_\sigma(t) = \sum_{k=0}^{\underline{\alpha}} \frac{(w^{(k)} * \phi_\sigma)(0)}{k!} t^k + I_{0+}^{\alpha+1/2} g_1 + I_{0+}^{\alpha+1/2} g_2.$$

Since $g_2$ is square integrable, the third term on the right belongs to the RKHS of $R^\alpha$, with squared RKHS-norm equal to $\|g_2\|_2^2 \leq \|I_{0+}^{3/2-\gamma}(w^{(\underline{\alpha})} * \phi_\sigma'')(t)\|_\infty^2$. This is $O(\sigma^{-1})$, by Lemma 5.2. The sum of the first two terms is a polynomial of degree $\underline{\alpha} + 1$ and the sum of squares of its coefficients is bounded by a constant times

$$\sum_{k=0}^{\underline{\alpha}} ((w^{(k)} * \phi_\sigma)(0))^2 + ((w^{(\underline{\alpha})} * \phi_\sigma')(0))^2.$$

The first term is bounded, while the second term is of order $\sigma^{2\gamma-2}$, which is $O(\sigma^{-1})$, since $\gamma > 1/2$.  $\square$

PROOF OF THEOREM 4.5.   The index $k$, when nested within a sum over $j$ in the following, is to be understood to range over all possible values $1, 2, \ldots, 2^{jd}$. The reproducing kernel Hilbert space of the variable $W$ is the set of functions $w = \sum_{j=1}^{J_\alpha} \sum_k w_{j,k} \psi_{j,k}$ with

$$\|w\|_{\mathbb{H}}^2 = \sum_{j=1}^{J_\alpha} \sum_k \frac{w_{j,k}^2}{\mu_j^2} < \infty.$$

For a fixed integer $J \leq J_\alpha$ (to be determined later), the projection $w_0^J = \sum_{j=1}^{J} \sum_k w_{0;j,k} \psi_{j,k}$ is clearly contained in the RKHS, whence, for any $\varepsilon > 0$,

$$\inf\{\|w\|_{\mathbb{H}}^2 : \|w - w_0^J\|_\infty < \varepsilon\}$$

$$\leq \|w_0^J\|_{\mathbb{H}}^2 = \sum_{j=1}^{J} \sum_k \frac{w_{0;j,k}^2}{\mu_j^2} \leq \sum_{j=1}^{J} 2^{j(2a-2\beta+d)} \|w_0\|_{\beta|\infty,\infty}^2$$

for $\mu_j 2^{jd/2} = 2^{-ja}$. For any numbers $\alpha_j \geq 0$ with $\sum_{j=1}^{J_\alpha} \alpha_j \leq 1$, we have

$$\Pr(\|W\|_\infty < \varepsilon) = \Pr\left(\sum_{j=1}^{J_\alpha} 2^{jd/2} \max_k |\mu_j Z_{j,k}| < \varepsilon\right)$$

$$\geq \prod_{j=1}^{J_\alpha} \prod_k \Pr(|\mu_j 2^{jd/2} Z_{j,k}| < \alpha_j \varepsilon).$$

Therefore, for $\mu_j 2^{jd/2} = 2^{-ja}$ and $\alpha_j = (K + d^2 j^2)^{-1}$ and a large constant $K$, it follows that

$$-\log \Pr(\|W\|_\infty < \varepsilon_n) \leq -\sum_{j=1}^{J_\alpha} 2^{jd} \log(2\Phi(\alpha_j \varepsilon_n 2^{ja}) - 1)$$

$$\lesssim \int_1^{2^{J_\alpha d}} -\log\left(2\Phi\left(\frac{\varepsilon_n x^{a/d}}{K + \log_2^2 x}\right) - 1\right) dx.$$

To justify the last step, we may choose the constant $K$ sufficiently large that the function $x \mapsto x^{a/d}/(K + \log_2^2 x)$ is nondecreasing on $[1, \infty)$.

The function $f : [0, \infty) \to \mathbb{R}$ defined by $f(y) = -\log(2\Phi(y) - 1)$ is decreasing from $\infty$ at $y = 0$ to $0$ at $y = \infty$. It is bounded above by a multiple of $1 + |\log y|$ for $y$ in an interval $[0, c]$ and bounded above by a multiple of $e^{-y^2/2}$ for $y \geq c$. [For the latter note that $f'(y) = -2\phi(y)/(2\Phi(y) - 1)$ is bounded above in absolute value by $2\phi(y)/(2\Phi(c) - 1)$ for $y \geq c$ so that $f(y) = f(\infty) - \int_y^\infty f'(x) \, dx$ is bounded in absolute value by $2(1 - \Phi(y))/(2\Phi(c) - 1)$ on this interval.]

We consider two cases to further bound the integral in the last display. For $\varepsilon_n 2^{J_\alpha a} \leq (K + J_\alpha^2 d^2)$, the argument $(\varepsilon_n x^{a/d})/(K + \log_2^2 x)$ is bounded above by a constant on the integration interval $[1, 2^{J_\alpha d}]$ and hence the function $f$ in the integral can be bounded above by a multiple of $1 + |\log|$, yielding as upper bound a multiple of

$$\int_1^{2^{J_\alpha d}} \left(1 + |\log|\left(\frac{\varepsilon_n x^{a/d}}{K + \log_2^2 x}\right)\right) dx \lesssim 2^{J_\alpha d}(\log(1/\varepsilon_n) + J_\alpha).$$

Whenever $a > 0$ and, in particular, if $\varepsilon_n 2^{J_\alpha a} > (K + J_\alpha^2 d^2)$, we can change variables $\varepsilon_n x^{a/d} = y$ and rewrite the integral as

$$\left(\frac{1}{\varepsilon_n}\right)^{d/a} \int_{\varepsilon_n}^{\varepsilon_n 2^{J_\alpha a}} f\left(\frac{y}{K + (d/a)^2(\log_2 y + \log_2(1/\varepsilon_n))^2}\right) \frac{d}{a} y^{d/a-1} \, dy.$$

The integral in this expression is bounded above by

$$\left[\int_0^{1/\varepsilon_n} f\left(\frac{y}{K + (2d/a)^2 \log_2^2(1/\varepsilon_n)}\right)\right.$$

$$\left. + \int_{1/\varepsilon_n}^\infty f\left(\frac{y}{K + (2d/a)^2 \log_2^2 y}\right)\right] \frac{d}{a} y^{d/a-1} \, dy$$

$$\leq \mu_n^{d/a} \int_0^{1/(\varepsilon_n \mu_n)} f(x) \frac{d}{a} x^{d/a-1} \, dx$$

$$+ \int_0^\infty f\left(\frac{y}{K + (2d/a)^2 \log_2^2 y}\right) \frac{d}{a} y^{d/a-1} \, dy$$

for $\mu_n = K + (2d/a)^2 \log_2^2(1/\varepsilon_n)$. The integral in the first term on the right is bounded as $n \to \infty$, whence the whole expression is bounded by a multiple of $(\log(1/\varepsilon_n))^{2d/a}$.

Combining the preceding, we conclude that

$$\phi_{w_0^J}(\varepsilon_n) = \inf\{\|w\|_{\mathbb{H}}^2 : \|w - w_0^J\|_\infty < \varepsilon_n\} - \log \Pr(\|W\|_\infty < \varepsilon_n)$$

$$\lesssim \sum_{j=1}^J 2^{j(2a-2\beta+d)} + \begin{cases} 2^{J_\alpha d}(\log(1/\varepsilon_n) + J_\alpha), & \text{if } \varepsilon_n 2^{J_\alpha a} \lesssim J_\alpha^2, \\ \left(\dfrac{1}{\varepsilon_n}\right)^{d/a}(\log(1/\varepsilon_n))^{2d/a}, & \text{if } \varepsilon_n 2^{J_\alpha a} \gtrsim J_\alpha^2. \end{cases}$$

The display gives the concentration function at the projection $w_0^J$. By Theorem 2.1, there exist measurable sets $B_n$ satisfying the three assertions of Theorem 4.5, but with $w_0$ replaced by $w_0^J$ and the 4 in the last condition replaced by 2. Since

$$\|w_0 - w_0^J\|_\infty \le \sum_{j=J+1}^\infty 2^{jd/2} \max_k |w_{0;j,k}|$$

$$\le \sum_{j=J+1}^\infty 2^{-j\beta} \|w_0\|_{\beta|\infty,\infty} \lesssim 2^{-J\beta},$$

we have the three assertions of Theorem 4.5 as given as soon as

$$\phi_{w_0^J}(\varepsilon_n) \le n\varepsilon_n^2 \quad \text{and} \quad 2^{-J\beta} \le \varepsilon_n.$$

The proof is completed by verifying that $\varepsilon_n$ as given in the theorem satisfies these inequalities in the various cases, for suitable $J$ (set $J = J_\alpha$ if $a \le \alpha$ and $J = J_a$ otherwise). We omit the (tedious) derivation of this. □

DEPARTMENT OF MATHEMATICS
VRIJE UNIVERSITEIT
DE BOELELAAN 1081A
1081 HV AMSTERDAM
THE NETHERLANDS
E-MAIL: aad@cs.vu.nl
          harry@cs.vu.nl